\documentclass[12pt]{amsproc}
\usepackage[cp1251]{inputenc}
\usepackage[T2A]{fontenc} 
\usepackage[russian]{babel} 
\usepackage{amsbsy,amsmath,amssymb,amsthm,graphicx,color} 
\usepackage{amscd}
\def\le{\leqslant}
\def\ge{\geqslant}

\newtheorem{prop}{Предложение}

\theoremstyle{definition}

\theoremstyle{remark}

\begin {document}
\unitlength=1mm
\title[Распределение по высоте в системах корней]
{Распределение по высоте в системах корней}
\author{Г. Г. Ильюта}
\email{ilgena@rambler.ru}
\address{}
\thanks{Работа поддержана грантом РФФИ-16-01-00409}

\bigskip

\begin{abstract}
Мы докажем несколько формул для распределения положительных корней.

We prove several formulas for the distribution of positive roots.
\end{abstract}

\maketitle
\tableofcontents

\bigskip

\section{Введение}

\bigskip

  Мы изучаем распределение по высоте корня в системах корней $A_n$, $B_n$, $C_n$, $D_n$, $E_6$, $E_7$, $E_8$, $F_4$, $G_2$. В любой из этих систем корней $\Phi$ выбор базиса из простых корней определяет на $\Phi$ функцию высоты: $ht(\alpha)$ для $\alpha\in\Phi$ равно сумме координат в разложении $\alpha$ по базису. Обозначим через $b_k$ количество корней высоты $k$ в системе корней $\Phi$. Поскольку $b_{-k}=b_k$ для всех $k$, то достаточно рассматривать только $k>0$. Мы получим несколько формул для (производящих функций) чисел $b_k$, используя хорошо известное равенство, связывающее эти числа и показатели Кокстера $e_1,\dots, e_n$ системы корней $\Phi$ (образованное показателями Кокстера разбиение сопряжено образованному числами $b_k$ разбиению, тем самым списки чисел $b_k$ легко восстанавливаются по приведённым в п. 2 спискам показателей Кокстера) \cite{11}
$$
B(q):=\sum_{\alpha\in\Phi^+}q^{ht(\alpha)-1}=\sum_{k=1}^{h-1}b_kq^{k-1}
$$
$$
=\frac{E(q)-n}{q-1}=\sum_{i=1}^n(1+q+\dots+q^{e_i-1}),    \eqno (1)
$$
где $\Phi^+$ -- множество положительных корней по отношению к выбранному базису, $h$ -- число Кокстера системы корней $\Phi$ (порядок элемента Кокстера), $E(q):=\sum_{i=1}^nq^{e_i}$. 

  В п. 2 числа $b_k$ представлены как скалярные произведения характеров группы Вейля. Это представление указывает на то, что характеры 
$$
R-\sum_{j=0}^{k-1}\chi_j\uparrow W,\qquad k=1,\dots,h,
$$
можно рассматривать как обобщения чисел $b_k$. Здесь для $j=0,1,\dots,h-1$ $\chi_j$ -- характер порождённой элементом Кокстера $C_{\Phi}$ циклической группы, $\chi_j(C_{\Phi}):=\epsilon^j$, $\chi_j\uparrow W$ -- индуцированный характер, $R$ -- регулярный характер группы Вейля, $\epsilon:=e^{\frac{2\pi i}{h}}$. Возникает вопрос, какие ещё связанные с системой корней конструкции и факты можно поднять в кольцо характеров группы Вейля. В п. 2 мы сопоставим характер группы Вейля любому симметрическому многочлену от корней характеристического многочлена $C_{\Phi}(q)$ элемента Кокстера, в частности, сопоставим $q$-характер группы Вейля самому многочлену $C_{\Phi}(q)$. Прообразы при этом соответствии равны скалярным произведениям образов и характера $\chi$ представления группы Вейля как группы, порождённой отражениями.

  Равенство в п. 2 получено как следствие формулы (1) и представления для $E(q)$ в виде скалярного произведений характеров группы Вейля \cite{25}. Скалярное произведение в пространстве характеров определяется как сумма по группе Вейля. Известна ещё одна формула для $E(q)$, в которой присутствует явно суммирование по группе Вейля $W$,
$$
E(q)=\sum_{w\in W}(-1)^{l(w)}P_q(w(\theta +\rho)-\rho),   
$$
где $l(w)$ -- длина $w$ (количество базисных отражений в приведённом разложении элемента $w\in W$), $P_q$ -- $q$-аналог функции разбиений Костанта, $\theta$ -- единственный корень максимальной высоты $h-1$, $\rho$ -- вектор Вейля (полусумма положительных корней). Эта формула связывает $E(q)$ с многочленами Костки-Фолкеса и многочленами Каждана-Люстига, ссылки имеются в \cite{19}. Интересно было бы найти прямую связь между этими двумя появлениями суммирования по группе Вейля в формулах для $E(q)$.

  В п. 3 изучается интерполяция на множестве корней степени $h$ из единицы и на множестве примитивных корней степени $h$ из единицы. В частности, в качестве интерполируемой функции рассматривается многочлен, коэффициенты которого определяют арифметическую функцию Коэна (пример: $E(q)$). По определению \cite{6}, \cite{13} для фиксированного $h\in\mathbb Z_{>0}$ такая функция $a(k)$ зависит только от наибольшего общего делителя $(k,h)$
$$
a(k)=a((k,h)), k\in\mathbb Z_{>0}.
$$
В п. 4-5 формулы из п. 3 применяются к двум связанным дискретным преобразованием Фурье арифметическим функциям: $m(k)$ -- кратность корня $\epsilon^k$ многочлена $C_{\Phi}(q)$, $p(k)$ -- сумма $k$-х степеней корней многочлена $C_{\Phi}(q)$. Эти функции являются функциями Коэна, что следует из формул (12) и (13).

  Корнями многочлена $C_{\Phi}(q)$ являются числа \cite{11}
$$
\epsilon^{e_1},\dots, \epsilon^{e_n}         
$$
и поэтому
$$
E(q)=\sum_{i=1}^nq^{e_i}=\sum_{k=0}^{h-1}m(k)q^k,        \eqno (2)
$$
Целочисленный многочлен, все корни которого являются корнями из единицы степени $h$, представим в следующем виде
$$
C_{\Phi}(q)=\prod_{d|h}(q^d-1)^{e(d)}.                    \eqno (3)
$$
В п. 4 приводится список таких представлений для всех систем корней и изучаются связи многочлена $B(q)$ с арифметическими функциями $e(d)$, $m(k)$ и $p(k)$. В случае $A-D-E$ геометрическое определение чисел $e(d)$ имеется в \cite{14}.

  Для $d,h\in\mathbb Z_{>0}$ определим функцию $a_{d,h}(x)$, $x\in\mathbb R_{>0}$, условием: $a_{d,h}(x)$ равно количеству целых чисел $j$, $1\le j\le x$, для которых $d|(j,h)$ -- наибольший общий делитель $j$ и $h$. В предельном случае $h=0$ функция $a_{d,h}(x)$ равна $[x/d]$ (целая часть $x/d$) и последний шаг в доказательстве Предложения 2 становится применением классического метода гиперболы Дирихле -- равенства, представляющего два способа перебора целых точек с положительными координатами под гиперболой $dd'=x$, 
$$
\sum_{n\le x}\sum_{d|n}f(d)=\sum_{m\le x}f(m)\left[\frac{x}{m}\right], 
$$
где $f(m)$ -- арифметическая функция. Традиционно этот метод использовался в аналитической теории чисел для асимптотических оценок \cite{4}, p. 65. Аналогичное нашему применение метода гиперболы имеется в \cite{24}, где доказано похожее на формулу (14) равенство для количества целых точек в некоторых эллипсах, связанных с мнимыми квадратичными полями.

  В п. 5 изучаются связи многочлена $B(q)$ с разложением Мунаги \cite{16}: для любого комплексного многочлена $\sum a(i)q^i$ степени $<h$ существует единственное представление вида
$$
\frac{a(0)+a(1)q+\dots+a(h-1)q^{h-1}}{1-q^h}=\sum_{d|h} \frac{H_d(q)}{1-q^d},                                                                   \eqno (4)
$$
где $\deg H_d(q)<\phi(d)$ для всех $d|h$, $\phi(d)$ -- функция Эйлера (количество примитивных корней степени $d$ из единицы). Мы докажем, что арифметическим функциям Коэна отвечают простейшие разложения Мунаги. Заметим также, что разложение Мунаги можно считать $q$-аналогом тождества Гаусса для функции Эйлера $h=\sum_{d|h}\phi(d)$. Вопрос: можно ли определить в терминах систем корней коэффициенты (многочлены $E_d(q)$ и $F_d(q)$, $d|h$) разложений Мунаги для многочленов $B(q)$ и $qB(q)$.

  Для случая $A-D-E$ в п. 6 содержится формула для многочлена $B(q)$, являющаяся следствием формулы (1) и представления для $E(q)$, зависящего от $h$ и от показателей квазиоднородности отвечающей системе корней простой особенности \cite{23}. Различные соотношения для показателей квазиоднородности позволяют представить эту формулу как зависящую от других параметров: порядка или размерности специального представления отвечающей системе корней бинарной полиэдральной группы, сферического объёма отвечающего системе корней треугольника Шварца, определителя матрицы Картана, длин ветвей диаграммы Дынкина.

  Для случая $A-D-E$ в п. 7 содержатся следствия формулы (1), связанные с многочленом Дынкина представления отвечающей системе корней группы Ли \cite{19}, p. 301, и перечислением антицепей для частичного порядка на $\Phi^+$ \cite{18}, р. 1206. 

  В п. 2-7 при поиске формул для $B(q)$ мы несколькими способами преобразовали правую часть формулы (1). Возможен другой подход к поиску содержащих числа $b_k$ производящих функций -- можно сделать в формуле (1) "замену базиса" $q^i\to f_i$, $i=0,1,\dots$. Мы ограничимся двумя примерами. Замена в формуле (1)
$$
q^i\to\log\frac{1-q^{i+2}}{1-q^{i+1}},\qquad i=0,1,\dots   
$$
приводит к хорошо известному равенству двух мультипликативных представлений для функции роста группы Вейля $W$ по отношению к множеству базисных отражений (или, в другой интерпретации, для многочлена Пуанкаре многообразия флагов)
$$
\sum_{w\in W}q^{l(w)}=\prod_{\alpha\in\Phi^+}\frac{1-q^{ht(\alpha)+1}}{1-q^{ht(\alpha)}}=\prod_{i=1}^n\frac{1-q^{e_i+1}}{1-q},                              \eqno (5)
$$
Следующие многочлены известны в теории чисел как многочлены Мириманова \cite{3}
$$
T_{k,m}(q):=\sum_{i=1}^ki^mq^i.
$$
В формулах для многочленов Мириманова появляются многочлены Апостола-Бернулли \cite{3}, многочлены Эйлера \cite{5}, p. 328, и числа Стирлинга \cite{5}, р. 326. Замена  $q^i\to i^mq^i$, $i=1,2,\dots$, в умноженной на $q$ формуле (1) приводит к равенству
$$
\left(q\frac{d}{dq}\right)^m(qB(q))=\sum_{k=1}^{h-1}b_kk^mq^k=\sum_{i=1}^nT_{e_i,m}(q).
$$
В качестве более общей замены можно использовать, например, $q^i\to f(i)q^i$, где $f(t)$ -- определённое число раз дифференцируемая функция. В \cite{10}, p. 30, имеются и цитируются формулы для сумм $\sum_if(i)q^i$, которые появляются в правой части формулы (1) при такой замене.

  Известны обобщения утверждения о сопряжённости разбиений, состоящих из показателей Кокстера и из чисел $b_k$ \cite{1}, \cite{2}, \cite{12}, \cite{19} (там же можно найти историю появления формул (1) и (5) в работах А. Шапиро, Р. Стейнберга, Б. Костанта и И. Макдональда), но мы ограничимся рассмотрением только классического случая. Конечные системы корней отвечают простым особенностям голоморфных функций. Для любых изолированных особенностей гиперповерхностей обобщением набора показателей Кокстера $e_1,\dots, e_n$ является подходящим образом нормализованный спектр особенности $s_1,\dots, s_n$ \cite{20}. Вопрос: существует ли для любой особенности интересная интерпретация коэффициентов двух обобщений многочлена $B(q)$ (собственные числа оператора монодромии особенности являются корнями из единицы степени $h$):
$$
\sum_{i=1}^n\frac{1-q^{s_i}}{1-q},\qquad
\sum_{i=1}^n\frac{1-q^{s_i\: mod\: h}}{1-q},
$$
где $0\le s_i\: mod\: h<h,\; i=1,\dots,n$.

\bigskip

\section{Индуцированные характеры группы Вейля}

\bigskip

  Пусть $(.,.)$ -- скалярное произведение в пространстве характеров группы Вейля (неприводимые характеры образуют ортонормированный базис).

\begin{prop}\label{prop1} Для $k=1,\dots,h$
$$
b_k=(R-\sum_{j=0}^{k-1}\chi_j\uparrow W,\chi)=(R-\sum_{j=1}^{k-1}\chi_j\uparrow W,\chi).                                                     \eqno (6)
$$
\end{prop}

  Доказательство. Согласно \cite{25} 
$$
\sum_{i=1}^nq^{e_i}=\sum_{k=0}^{h-1}(\chi_k\uparrow W,\chi)q^k,
$$
что эквивалентно равенству
$$
m(k)=(\chi_k\uparrow W,\chi),\qquad k=0,\dots, h-1.
$$
Скалярное произведение регулярного и любого неприводимого характеров равно размерности этого неприводимого характера и $\chi$ является неприводимым характером размерности $n$. Поэтому
$$
(R-\sum_{j=0}^{k-1}\chi_j\uparrow W,\chi)=n-\sum_{j=0}^{k-1}m(j)=n-\sum_{j=1}^{k-1}m(j)=b_k.
$$
Последнее равенство доказано в Предложении 8.$\blacksquare$

  Определим циклотомический многочлен $\Phi_h(q)$ и суммы Рамануджана $c_h(j)$ -- степенные суммы корней многочлена $\Phi_h(q)$ \cite{13}. 
$$
\Phi_h(q):=\prod_{i=1}^{\phi(h)}(q-\epsilon^{k_i})
=\prod_{d|h}(q^d-1)^{\mu\left(\frac{h}{d}\right)},          
$$
$$
q^h-1=\prod_{d|h}\Phi_d(q),          
$$
где $1\le k_1,\dots,k_{\phi(h)}\le h-1$ -- взаимно простые с $h$ числа, $\mu(k)$ -- функция Мёбиуса.
$$
c_h(j):=\sum_{i=1}^{\phi(h)}\epsilon^{jk_i}
=\sum_{d|(j,h)}d\mu\left(\frac{h}{d}\right)
=\frac{\phi(h)\mu\left(\frac{h}{(j,h)}\right)}{\phi\left(\frac{h}{(j,h)}\right)}.                                       \eqno (7)
$$
В частности, для функции Эйлера $\phi(h)$ и функции Мёбиуса $\mu(h)$ имеем
$$
\phi(h)=c_h(0)=\sum_{i=1}^{\phi(h)}1=\sum_{d|h}d\mu\left(\frac{h}{d}\right),
$$
$$
\mu(h)=c_h(1)=\sum_{i=1}^{\phi(h)}\epsilon^{k_i}.
$$

  Формула (6) показывает, что характеры $R-\sum_{j=0}^{k-1}\chi_j\uparrow W$
можно рассматривать как обобщения чисел $b_k$. Кратностям $m(k)$ в этом смысле отвечают характеры $\chi_k\uparrow W$. По определению для $k\in\mathbb Z$
$$
p(k)=\sum_{j=0}^{h-1}m(j)\epsilon^{jk}
$$
$$
=\sum_{d|h}m(d)\sum_{(j,h)=d}\epsilon^{jk}=\sum_{d|h}m(d)\sum_{(j/d,h/d)=1}\epsilon^{jk}
$$
$$
=\sum_{d|h}m(d)c_{h/d}(k)=\sum_{d|h}m(h/d)c_d(k),
$$
другими словами, функция $p(k)$ является дискретным преобразованием Фурье функции $m(k)$. Поэтому степенным суммам $p(k)$ можно сопоставить характеры
$$
\sum_{d|h}c_d(k)\chi_{h/d}\uparrow W.              
$$
Обратное преобразование Фурье имеет вид
$$
m(k)=\frac{1}{h}\sum_{j=0}^{h-1}p(j)\epsilon^{-jk}
$$
$$
=\frac{1}{h}\sum_{d|h}p(d)c_{h/d}(-k)=\frac{1}{h}\sum_{d|h}p(h/d)c_d(k).
$$

  Степенные суммы порождают кольцо симметрических функций. Пусть произведению в этом кольце отвечает произведение характеров (тензорное произведение представлений). Тогда можно использовать характеры, отвечающие степенным суммам $p(k)$, чтобы сопоставить характер любому симметрическому многочлену от корней многочлена $C_{\Phi}(q)$, например, коэффициентам этого многочлена. Тем самым $q$-характер группы Вейля будет сопоставлен и самому многочлену $C_{\Phi}(q)$. Известны формулы, связывающие коэффициенты многочлена и степенные суммы его корней. В Предложении 2 содержатся равенства, связывающие коэффициенты многочлена $C_{\Phi}(q)$ с арифметическими функциями $e(d)$ и $m(k)$.

  Следующий многочлен является $q$-аналогом функции Эйлера
$$
\Psi_h(q):=\sum_{i=1}^{\phi(h)}q^{k_i}.
$$
В Предложении 6 нам понадобится равенство
$$
\Psi_h(\epsilon^k)=c_h(k),\quad k\in\mathbb Z.
$$

\begin{prop}\label{prop2}
$$
\frac{C'_{\Phi}(q)}{C_{\Phi}(q)}=\sum_{d|h}e(d)\frac{dq^{d-1}}{q^d-1}  
$$
$$
=\sum_{k=0}^{h-1}\frac{m(k)}{q-\epsilon^k}  
$$
$$
=\sum_{d|h}m\left(\frac{h}{d}\right)\frac{\Phi'_d(q)}{\Phi_d(q)}  
$$
$$
=\sum_{d|h}m\left(\frac{h}{d}\right)\sum_{d'|d}\mu\left(\frac{d}{d'}\right)\frac{d'q^{d'-1}}{q^{d'}-1}
$$
$$
=\sum_{d|h}m\left(\frac{h}{d}\right)\frac{\sum_{j=1}^{d}c_d(j)q^{j-1}}{q^d-1}  
$$
$$
=\sum_{d|h}m\left(\frac{h}{d}\right)\frac{\phi(d)}{q^d-1}\sum_{d'|d}\frac{\mu(d')}{\phi(d')}\frac{\Psi_{d'}(q^{d/d'})}{q}.  
$$
\end{prop}

  Доказательство. Формулу (3) можно переписать в следующем виде
$$
C_{\Phi}(q)=\prod_{k=0}^{h-1}(q-\epsilon^k)^{m(k)}=\prod_{d|h}\Phi_d^{m(\frac{h}{d})}(q).           
$$
Применяя логарифмическую производную к этим равенствам, получим первые три формулы из Предложения 2. Остаётся применить равенства для $\Phi'_d(q)/\Phi_d(q)$ из \cite{17}, \cite{15} и \cite{26}
$$
\frac{\Phi'_d(q)}{\Phi_d(q)}=\sum_{d'|d}\mu\left(\frac{d}{d'}\right)\frac{d'q^{d'-1}}{q^{d'}-1}  
$$
$$
=\frac{\sum_{j=1}^{d}c_d(j)q^{j-1}}{q^d-1}  
$$
$$
=\frac{\phi(d)}{q^d-1}\sum_{d'|d}\frac{\mu(d')}{\phi(d')}\frac{\Psi_{d'}(q^{d/d'})}{q}.\blacksquare  
$$

\bigskip

\section{Интерполяция в корнях из единицы}

\bigskip

  Обозначим через $L[F(q)](q)$ ($L^*[F(q)](q)$) интерполяционный многочлен Лагранжа функции $F(q)$ в корнях степени $h$ из единицы (примитивных корнях степени $h$ из единицы),
$$
L[F(q)](q):=(q^h-1)\frac{1}{h}\sum_{i=0}^{h-1}\frac{\epsilon^iF(\epsilon^i)}{(q-\epsilon^i)},                    \eqno (8)
$$
$$
L^*[F(q)](q):=\Phi_h(q)\sum_{i=1}^{\phi(h)}\frac{F(\epsilon^{k_i})}{\Phi'_h(\epsilon^{k_i})(q-\epsilon^{k_i})},  \eqno (9)
$$
Если $F(q)$ является многочленом, то ниже без дополнительных ссылок мы используем следующий факт: интерполяционный многочлен для $F(q)$ совпадает с $F(q)$, если $\deg F(q)$ меньше количества узлов интерполяции.

Согласно \cite{7}
$$
L[F(q)](q)=\frac{1}{h}\sum_{i=0}^{h-1}F(\epsilon^i)\sum_{k=0}^{h-1}\epsilon^{(h-i)k}q^k
$$

  Пусть 
$$
D:=diag(h,\dots,h),\quad \dim D=h\times h,
$$
$$
F_h:=(F(1),F(\epsilon),\dots,F(\epsilon^{h-1})),
$$
$$
Q_k:=(1,q,\dots,q^{k-1}),\quad k\in\mathbb Z_{>0},
$$
$$
V=(\epsilon^{(i-1)(j-1)})_{i,j=1,\dots,h}.
$$

\begin{prop}\label{prop3} 
$$
L[F(q)](q)=-\frac{1}{h^h}\det\begin{pmatrix}0 & Q_h\\\overline V^tF_h^t & D\end{pmatrix}.  
$$
\end{prop}

  Доказательство. Известное детерминантное представление общего многочлена Лагранжа в нашем случае имеет вид
$$
L[F(q)](q)=-\frac{\det\begin{pmatrix}0 & Q_h\\F_h^t & V\end{pmatrix}}{\det V}. 
$$
Умножаем числитель дроби в правой части этой формулы на
$$
\det\begin{pmatrix}1 & 0\\0 & \overline V^t\end{pmatrix},  
$$
а знаменатель -- на $\det \overline V^t$. Из соотношений ортогональности для характеров циклической группы следует, что 
$$
\overline V^tV=D, \quad \det (\overline V^tV)=h^h.\blacksquare
$$

Пусть 
$$
C:=(c_h(i+j-2))_{i,j=1,\dots,\phi(h)},
$$
$$
F^*_h:=(F(\epsilon^{k_1}),\dots,F(\epsilon^{k_{\phi(h)}})),
$$
$$
V^*=(\epsilon^{k_i(j-1)})_{i,j=1,\dots,\phi(h)}.
$$

\begin{prop}\label{prop4} Для $h>1$
$$
L^*[F(q)](q)=\frac{\prod_{p|h,p\, prime}p^{\frac{\phi(h)}{p-1}}}{(-1)^{1+\phi(h)/2}h^{\phi(h)}}\det\begin{pmatrix}0 & Q_{\phi(h)}\\V^{*t}F_h^{*t} & C\end{pmatrix}. 
$$
\end{prop}

  Доказательство. В правой части формулы
$$
L^*[F(q)](q)=-\frac{\det\begin{pmatrix}0 & Q_{\phi(h)}\\F_h^{*t} & V^*\end{pmatrix}}{\det V^*}  
$$
умножаем числитель дроби на
$$
\det\begin{pmatrix}1 & 0\\0 & V^{*t}\end{pmatrix},  
$$
а знаменатель -- на $\det V^{*t}$. Используем равенства 
$$
V^{*t}V^*=C,\quad\det (V^{*t}V^*)=\det C=discr \Phi_h(q).
$$ 
Для дискриминанта многочлена $\Phi_h(q)$ имеем формулу \cite{21}
$$
discr \Phi_h(q)=\frac{(-1)^{\phi(h)/2}h^{\phi(h)}}{\prod_{p|h,p\, prime}p^{\frac{\phi(h)}{p-1}}}.\blacksquare 
$$

\begin{prop}\label{prop5} Если 
$$
A(q):=\sum_{k=0}^{h-1}a(k)q^k
$$
и $a(k)$ является арифметической функцией Коэна, то
$$
\frac{A(q)}{1-q^h}=\frac{1}{h}\sum_{d|h}A(\epsilon^{h/d})\frac{\sum_{k=0}^{d-1}c_d(k)q^k}{1-q^d} 
$$
$$
=\frac{1}{h}\sum_{d|h}A(\epsilon^{h/d})\frac{\Phi'_d(1/q)}{q\Phi_d(1/q)}      
$$
$$
=\frac{1}{h}\sum_{d|h}A(\epsilon^{h/d})\sum_{d'|d}\frac{d'\mu(d/d')}{1-q^{d'}}      
$$
$$
=\frac{1}{h}\sum_{d|h}\frac{A(\epsilon^{h/d})\phi(d)}{1-q^d} \left(1-q^d+\sum_{d'|d}\frac{\mu(d')}{\phi(d')}\Psi_{d'}(q^{d/d'})\right).      $$
\end{prop}

  Доказательство. В \cite{7} доказано следующее равенство
$$
A(q)=\frac{1}{h}\sum_{k=0}^{h-1}q^k\sum_{d|h}A(\epsilon^{h/d})c_d(k)            $$
Меняя порядок суммирования и используя $d$-периодичность сумм Рамануджана $c_d(k)$, имеем
$$
\frac{A(q)}{1-q^h}=\frac{1}{h}\frac{\sum_{d|h}A(\epsilon^{h/d})\sum_{k=0}^{h-1}c_d(k)q^k}{1-q^h}       
$$
$$
=\frac{1}{h}\frac{\sum_{d|h}A(\epsilon^{h/d})(1+q^d+\dots+q^{h/d-1})\sum_{k=0}^{d-1}c_d(k)q^k}{1-q^h}
$$
$$
=\frac{1}{h}\sum_{d|h}A(\epsilon^{h/d})\frac{\sum_{k=0}^{d-1}c_d(k)q^k}{1-q^d}.$$
Остаётся применить формулы для $\sum_{k=0}^{d-1}c_d(k)q^k$ из \cite{17}, \cite{15} и \cite{26}
$$
\frac{\sum_{k=0}^{d-1}c_d(k)q^k}{1-q^d}=\frac{\Phi'_d(1/q)}{q\Phi_d(1/q)}
$$
$$
=\sum_{d'|d}\frac{d'\mu(d/d')}{1-q^{d'}}
$$
$$
=\frac{\phi(d)}{1-q^d} \left(1-q^d+\sum_{d'|d}\frac{\mu(d')}{\phi(d')}\Psi_{d'}(q^{d/d'})\right).\blacksquare
$$

 В Предложении 6 мы специализируем для $A(q)=\Psi_h(q)$ формулы из \cite{7} и Предложения 5. В этом случае согласно формуле (7)
$$
A(\epsilon^{h/d})=\Psi_h(\epsilon^{h/d})=c_h(h/d)=\frac{\phi(h)\mu(d)}{\phi(d)}.
$$
По-видимому, некоторые равенства из Предложения 6 можно упростить, например, с помощью соотношений ортогональности для сумм Рамануджана. В \cite{27}, p. 75, доказана формула
$$
\frac{\Psi_h(q)}{1-q^h}=\sum_{d|h}\frac{\mu(d)q^{d}}{1-q^{d}}=\sum_{d|h}\frac{\mu(d)}{1-q^{d}}-\sum_{d|h}\mu(d),                                \eqno (10)
$$
для $h>1$
$$
\frac{\Psi_h(q)}{1-q^h}=\sum_{d|h}\frac{\mu(d)}{1-q^{d}}.        \eqno (11)
$$

\begin{prop}\label{prop6} 
$$
\frac{\Psi_h(q)}{1-q^h}=\frac{1}{h}\sum_{k=0}^{h-1}\frac{\epsilon^kc_h(k)}{\epsilon^k-q}
$$
$$
=\frac{1}{h}\frac{\sum_{i=0}^{h-1}c_h(i)\sum_{k=0}^{h-1}\epsilon^{(h-i)k}q^k}{1-q^h}
$$
$$
=\frac{1}{h}\frac{\sum_{k=0}^{h-1}q^k\sum_{d|h}c_h(h/d)c_d(k)}{1-q^h}           $$
$$
=\frac{1}{h}\sum_{d|h}c_h(h/d)\frac{\sum_{k=0}^{d-1}c_d(k)q^k}{1-q^d} 
$$
$$
=\frac{1}{h}\sum_{d|h}c_h(h/d)\frac{\Phi'_d(1/q)}{q\Phi_d(1/q)}      
$$
$$
=\frac{1}{h}\sum_{d|h}c_h(h/d)\sum_{d'|d}\frac{d'\mu(d/d')}{1-q^{d'}}      
$$
$$
=\frac{\phi(h)}{h}\sum_{d|h}\frac{\mu(d)}{1-q^d} \left(1-q^d+\sum_{d'|d}\frac{\mu(d')}{\phi(d')}\Psi_{d'}(q^{d/d'})\right).\blacksquare      
$$
\end{prop}

\begin{prop}\label{prop7} В условиях Предложения 5
$$
\frac{A(q)}{1-q^h}-a(0)=\frac{\sum_{d|h}a(h/d)\Psi_d(q^{h/d})}{1-q^h}
$$
$$
=\sum_{d|h}a(h/d)\sum_{d'|d}\mu(d')\frac{q^{hd'/d}}{1-q^{hd'/d}}
$$
$$
=a(0)\frac{q^h}{1-q^h}+\sum\limits_{\substack{d|h\\d>1}}a(h/d)\sum_{d'|d}\frac{\mu(d')}{1-q^{hd'/d}}.
$$
$$
\frac{A(q)}{1-q^h}=\sum_{d|h}a(h/d)\sum_{d'|d}\frac{\mu(d')}{1-q^{hd'/d}}.
$$
\end{prop}

Доказательство.
$$
A(q)-a(0)(1-q^h)=\sum_{k=1}^ha(k)q^k
$$
$$
=\sum_{d|h}a(d)\sum_{(k,h)=d}q^k=\sum_{d|h}a(d)\sum_{(k/d,h/d)=1}q^k
$$
$$
\sum_{d|h}a(d)\Psi_{h/d}(q^d)=\sum_{d|h}a(h/d)\Psi_d(q^{h/d})
$$
и используем формулы (10) и (11).$\blacksquare$

\bigskip

\section{Арифметические функции Коэна и элемент Кокстера}

\bigskip

  Приведём список показателей Кокстера и многочленов $C_{\Phi}(q)$, представленных в виде (3).

  $A_n$, $h=n+1$, показатели Кокстера $1, 2,\dots, n$,
$$
C_{A_n}(q)=\frac{q^{n+1}-1}{q-1}.
$$

  $B_n$, $C_n$, $h=2n$, показатели Кокстера $1, 3,\dots, 2n-1$,
$$
C_{B_n}(q)=C_{C_n}(q)=\frac{q^{2n}-1}{q^n-1}.
$$

  $D_n$, $h=2n-2$, показатели Кокстера $1, 3,\dots, 2n-3, n-1$,
$$
C_{D_n}(q)=\frac{(q^{2n-2}-1)(q^2-1)}{(q^{n-1}-1)(q-1)}.
$$

  $E_6$, $h=12$, показатели Кокстера $1, 4, 5, 7, 8, 11$,
$$
C_{E_6}(q)=\frac{(q^{12}-1)(q^3-1)(q^2-1)}{(q^6-1)(q^4-1)(q-1)}.
$$

  $E_7$, $h=18$, показатели Кокстера $1, 5, 7, 9, 11, 13, 17$,
$$
C_{E_7}(q)=\frac{(q^{18}-1)(q^3-1)(q^2-1)}{(q^9-1)(q^6-1)(q-1)}.
$$

  $E_8$, $h=30$, показатели Кокстера $1, 7, 11, 13, 17, 19, 23, 29$,
$$
C_{E_8}(q)=\frac{(q^{30}-1)(q^5-1)(q^3-1)(q^2-1)}{(q^{15}-1)(q^{10}-1)(q^6-1)(q-1)}.
$$

  $F_4$, $h=12$, показатели Кокстера $1, 5, 7, 11$,
$$
C_{F_4}(q)=\frac{(q^{12}-1)(q^2-1)}{(q^6-1)(q^4-1)}.
$$

  $G_2$, $h=6$, показатели Кокстера $1, 5$,
$$
C_{G_2}(q)=\frac{(q^6-1)(q-1)}{(q^3-1)(q^2-1)}.
$$

  Число $\epsilon^k$ является корнем многочлена $q^{h/d}-1$, $d|h$, если и только если $d|(k,h)$. Поэтому
$$
m(k)=\sum_{d|(k,h)}e\left(\frac{h}{d}\right),  k\in\mathbb Z.     \eqno (12)
$$
По определению
$$
p(k)=\sum_{d|h}e(d)\sum_{l=0}^{d-1}e^{\frac{2\pi ilk}{d}}, k\in\mathbb Z,
$$
и поэтому из соотношений ортогональности для характеров циклической группы вытекает равенство
$$
p(k)=\sum_{d|(k,h)}de(d),  k\in\mathbb Z.                          \eqno (13)
$$

\begin{prop}\label{prop8} Для $k=1,\dots,h-1$
$$
b_k=n-\sum_{i=0}^{k-1}m(i)=n-\sum_{d|h}a_{d,h}(k-1)e\left(\frac{h}{d}\right).                                    \eqno (14)          
$$
\end{prop}

Доказательство. Для любой системы корней минимальный показатель равен $1$, а максимальный -- $h-1$. Поэтому 
$$
m(0)=\sum_{d|h}e(d)=0,
$$  
$$
\sum_{k=0}^{h-1}m(k)=\deg C_{\Phi}(q)=n=p(0)=\sum_{d|h}de(d).
$$
Формулу (1) можно переписать в следующем виде
$$
B(q)=\sum_{k=1}^{h-1}m(k)\frac{1-q^k}{1-q}
$$
$$
=\sum_{k=1}^{h-1}\left(\sum_{j=k}^{h-1}m(j)\right)q^{k-1}=\sum_{k=1}^{h-1}\left(\sum_{j=0}^{h-1}m(j)-\sum_{j=0}^{k-1}m(j)\right)q^{k-1}
$$
$$
=\sum_{k=1}^{h-1}\left(n-\sum_{j=0}^{k-1}m(j)\right)q^{k-1}.
$$
Для $0<k<h$ имеем
$$
\sum_{j=0}^km(j)=\sum_{j=1}^k\sum_{d|(j,h)}e\left(\frac{h}{d}\right)=\sum_{d|h}a_{d,h}(k)e\left(\frac{h}{d}\right).\blacksquare                       
$$

  Полагая в Предложении 5 $A(q)=n-E(q)=p(0)-E(q)$, получим равенство
$$
A(\epsilon^k)=p(0)-p(k),\quad k\in\mathbb Z.
$$
Поэтому из формулы (1), \cite{7} и Предложения 5 вытекает

\begin{prop}\label{prop9}
$$
\frac{1-q}{1-q^h}B(q)=\frac{(n-E(q))}{1-q^h}
$$
$$
=\frac{1}{h}\sum_{k=1}^{h-1}\epsilon^k\frac{p(0)-p(k)}{\epsilon^k-q}  
$$
$$
=\frac{1}{h}\frac{\sum_{i=0}^{h-1}(p(0)-p(i))\sum_{k=0}^{h-1}\epsilon^{(h-i)k}q^k}{1-q^h}
$$
$$
=\frac{1}{h}\frac{\sum_{k=0}^{h-1}q^k\sum_{d|h}(p(0)-p(h/d))c_d(k)}{1-q^h}      $$
$$
=\frac{1}{h}\sum_{d|h}(p(0)-p(h/d))\frac{\sum_{k=0}^{d-1}c_d(k)q^k}{1-q^d} 
$$
$$
=\frac{1}{h}\sum_{d|h}(p(0)-p(h/d))\frac{\Phi'_d(1/q)}{q\Phi_d(1/q)}      
$$
$$
=\frac{1}{h}\sum_{d|h}(p(0)-p(h/d))\sum_{d'|d}\frac{d'\mu(d/d')}{1-q^{d'}}      
$$
$$
=\frac{1}{h}\sum_{d|h}\frac{(p(0)-p(h/d))\phi(d)}{1-q^d} \left(1-q^d+\sum_{d'|d}\frac{\mu(d')}{\phi(d')}\Psi_{d'}(q^{d/d'})\right).\blacksquare      
$$
\end{prop}

  Из формулы (1) и Предложения 7 для $A(q)=E(q)$ вытекает Предложение 10. Используем равенства $m(0)=0$,
$$
n=E(1)=\sum_{d|h}a(h/d)\Psi_d(1),
$$ 
$$
\Psi_d(1)=\phi(d)=\sum_{d'|d}\mu(d')\frac{d}{d'}.
$$

\begin{prop}\label{prop10}
$$
B(q)=\frac{n-E(q)}{1-q}
$$
$$
=\sum\limits_{\substack{d|h\\d>1}}m(h/d)\frac{\Psi_d(1)-\Psi_d(q^{h/d})}{1-q}
$$
$$
=\sum\limits_{\substack{d|h\\d>1}}m(h/d)\sum_{d'|d}\frac{\mu(d')}{1-q}\left(\frac{d}{d'}-\frac{1-q^h}{1-q^{hd'/d}}\right)
$$
$$
=\sum\limits_{\substack{d|h\\d>1}}m(h/d)\sum_{d'|d}\frac{\mu(d')q^{hd'/d}}{1-q}\left(\frac{d}{d'}-\frac{1-q^h}{1-q^{hd'/d}}\right).\blacksquare
$$
\end{prop}

\bigskip

\section{Разложение на $q$-простейшие дроби Мунаги}

\bigskip

\begin{prop}\label{prop11} Функция $a(i)$ является функцией Коэна, если и только если в разложении (4) $\deg H_d(q)=0$ для всех $d|h$, в частности,
$$
\frac{m(0)+m(1)q+\dots+m(h-1)q^{h-1}}{1-q^h}=\sum_{d|h} \frac{e(h/d)}{1-q^d}, \eqno (15)
$$        
$$
\frac{p(0)+p(1)q+\dots+p(h-1)q^{h-1}}{1-q^h}=\sum_{d|h} \frac{de(d)}{1-q^d}. \eqno (16)
$$
\end{prop}     

  Доказательство. Если $\deg H_d(q)=0$ для всех $d|h$, то, умножая формулу (4) на $1-q^h$ и раскрывая скобки, получим, что $a(k)=\sum_{d|(k,h)}H_d(q)$, и значит $a(k)$ является функцией Коэна. Обратно, если $a(k)$ является функцией Коэна, то существует функция $l(d)$, для которой $a(k)=\sum_{d|(k,h)}l(d)$ \cite{6}, \cite{13}. Полагая $H_d(q)=l(d)$, приходим к разложению вида (4). Формулы (15) и (16) вытекают из (12) и (13). Разложение Мунаги для функции Коэна $a(k)$ также легко получить, используя известное представление общего преобразования Мёбиуса $f(k)=\sum_{d|k}g(d)$, $k\in\mathbb Z_{>0}$, как переписывания ряда Ламберта в виде степенного ряда:
$$
\sum_{k=1}^\infty f(k)q^k=\sum_{m=1}^\infty \frac{g(m)q^m}{1-q^m}. 
$$
Поэтому
$$
\sum_{k=1}^\infty a(k)q^k=\sum_{d|h} \frac{l(d)q^d}{1-q^d}
$$
$$
=\sum_{d|h} \frac{l(d)}{1-q^d}-\sum_{d|h}l(d)=\sum_{d|h} \frac{l(d)}{1-q^d}-a(0).
$$
Из $h$-периодичности последовательности $a(k)$ следует, что
$$
\sum_{k=1}^\infty a(k)q^k=\frac{\sum_{k=1}^h a(k)q^k}{1-q^h} 
$$
$$
=\frac{\sum_{k=0}^{h-1} a(k)q^k}{1-q^h}-\frac{a(0)-a(h)q^h}{1-q^h}=\frac{\sum_{k=0}^{h-1} a(k)q^k}{1-q^h}-a(0) 
$$
Сравнивая последние два соотношения, получим
$$
\frac{\sum_{k=0}^{h-1} a(k)q^k}{1-q^h}=\sum_{d|h}\frac{l(d)}{1-q^d}.\blacksquare 
$$

  Из формул (1), (2) и (15) вытекает

\begin{prop}\label{prop12}
$$
B(q)=\frac{n}{1-q}-\frac{1-q^h}{1-q}\sum_{d|h} \frac{e(h/d)}{1-q^d}.\blacksquare               \eqno (17) 
$$
\end{prop}

\begin{prop}\label{prop13} Пусть для многочлена $B(q)$ разложение Мунаги (4) имеет вид
$$
\frac{b_1+b_2q+\dots+b_{h-1}q^{h-2}}{1-q^h}=\sum_{d|h} \frac{E_d(q)}{1-q^d} \eqno (18) 
$$
и $E^*_d$ -- коэффициент при $q^{\phi(d)-1}$ в многочлене $E_d(q)$. Тогда
$$
n=\sum_{d|h} E_d(0),  
$$
$$
n-1=\sum_{d|h} E_d'(0),  
$$
$$
1=\sum\limits_{\substack{d|h\\d\,prime}}E^*_d,
$$
$$
0=E^*_1=E_1(q).
$$
\end{prop}

  Доказательство. Умножаем формулу (18) на $1-q^h$ и сравниваем коэффициенты при $q^0,q,q^{h-2},q^{h-1}$, используем известные формулы для чисел $b_k$: $b_1=n$, $b_2=n-1$, $b_{h-1}=1$, $b_h=0$. Например, вклад в коэффициент при $q^{h-2}$ в правой части могут дать только те слагаемые, для которых
$$
h-d+\phi(d)-1\ge h-2.
$$
Решениями этого неравенства являются только простые числа и $1$, но $E_1(q)=0$ и поэтому вклад в коэффициент дают только $E^*_d$ для простых $d|h$.$\blacksquare$

\begin{prop}\label{prop14} Для $h>1$
$$
\frac{E_h(q)}{\Phi_h(q)}=(n-e(1))\frac{L^*[(1-q)^{-1}](q)}{\Phi_h(q)}
$$
$$
=(n-e(1))\sum_{i=1}^{\phi(h)}\frac{1}{\Phi_h'(\epsilon^{k_i})(1-\epsilon^{k_i})(q-\epsilon^{k_i})} 
$$
$$
=\frac{(n-e(1))}{\Phi_h(q)}\frac{\prod_{d|h,d\, prime}p^{\frac{\phi(h)}{p-1}}}{(-1)^{1+\phi(h)/2}h^{\phi(h)}}\det\begin{pmatrix}0 & Q_{\phi(h)}\\L_h^t & C\end{pmatrix},
$$
где у вектора $L_h$ для $j=1,\dots,\phi(h)$ $j$-я координата $L_{h,j}$ имеет вид
$$
L_{h,j}=\sum_{i=1}^{\phi(h)}\frac{\epsilon^{k_i(j-1)}}{1-\epsilon^{k_i}}=\frac{L^*[q^{j-1}\Phi_h'(q)](1)}{\Phi_h(1)}.
$$
\end{prop}

  Доказательство. Для $i=1,\dots,\phi(h)$ и $d<h$ число $\epsilon^{k_i}$ не является корнем многочлена $1-q^d$. Поэтому
$$
\left.\frac{1-q^h}{1-q^d}\right|_{q=\epsilon^{k_i}}=0.
$$
Тогда из формул (17) и (18) следуют равенства
$$
L^*[B(q)](q)=(n-e(1))L^*[(1-q)^{-1}](q),
$$
$$
L^*[B(q)](q)=L^*[E_h(q)](q)=E_h(q).
$$
Остаётся использовать формулу (9) и Предложение 4.$\blacksquare$

\begin{prop}\label{prop15} Для $1\le m\le h$
$$
\sum\limits_{\substack{d|h\\d>1}}L_{d,m+1}=m-\frac{h+1}{2}.
$$
\end{prop}

  Доказательство. 
$$
\sum_{k=1}^{h-1}\frac{\epsilon^{km}}{1-\epsilon^k}=\sum\limits_{\substack{d|h\\d<h}}\sum_{(k,h)=d}\frac{\epsilon^{km}}{1-\epsilon^k}=\sum\limits_{\substack{d|h\\d<h}}\sum_{(k/d,h/d)=1}\frac{\epsilon^{km}}{1-\epsilon^k}
$$
$$
=\sum\limits_{\substack{d|h\\d<h}}L_{h/d,m+1}=\sum\limits_{\substack{d|h\\d>1}}L_{d,m+1}.
$$
Согласно \cite{8}, p. 10, для $1\le m\le h$
$$
\sum_{k=1}^{h-1}\frac{\epsilon^{km}}{1-\epsilon^k}=m-\frac{h+1}{2}.\blacksquare
$$

\begin{prop}\label{prop16} 
$$
\frac{n}{1-q}=\sum_{d|h} \frac{E_d(q)+q^{d-1}E_d(1/q)}{1-q^d}.
$$
\end{prop}

  Доказательство. Известно, что соотношения для показателей Кокстера $e_k+e_{n+1-k}=h$ влекут соотношения $b_k+b_{h+1-k}=n$. А именно, из формулы (1) следует, что
$$
B(q)+q^{h-1}B(1/q)=\sum_{i=1}^n\left(\frac{1-q^{e_i}}{1-q}+q^{h-1}\frac{1-1/q^{h-e_i}}{1-1/q}\right)=n\frac{1-q^h}{1-q}.
$$
Используя формулу (18), получим равенство
$$
B(q)+q^{h-1}B(1/q)=\sum_{d|h}\left( \frac{1-q^h}{1-q^d}E_d(q)+q^{h-1}\frac{1-1/q^h}{1-1/q^d}E_d(1/q)\right)
$$
$$
=(1-q^h)\sum_{d|h} \frac{E_d(q)+q^{d-1}E_d(1/q)}{1-q^d}.\blacksquare
$$

  Предложения 17, 18 и 19 являются аналогами Предложений 13, 14 и 16 для многочлена $qB(q)$, доказательства также аналогичны.

\begin{prop}\label{prop17} Пусть для многочлена $qB(q)$ разложение Мунаги (4) имеет вид
$$
\frac{b_1q+b_2q^2+\dots+b_{h-1}q^{h-1}}{1-q^h}=\sum_{d|h} \frac{F_d(q)}{1-q^d}  
$$
и $F_2(q)=Aq+B$ (если $h$ нечётно, то $A=B=0$). Тогда
$$
n=1+\sum_{d|h} F'_d(0)=1+A+\sum_{d|h,d>2} F'_d(0),  
$$
$$
n-1=1+B+\frac{1}{2}\sum_{d|h,d>2} F''_d(0),  
$$
$$
1=F_1(q),  
$$
$$
0=\sum_{d|h} F_d(0).\blacksquare  
$$
\end{prop}

\begin{prop}\label{prop18} Для $h>1$
$$
\frac{F_h(q)}{\Phi_h(q)}=(n-e(1))\frac{L^*[q(1-q)^{-1}](q)}{\Phi_h(q)}
$$
$$
=(n-e(1))\sum_{i=1}^{\phi(h)}\frac{\epsilon^{k_i}}{\Phi_h'(\epsilon^{k_i})(1-\epsilon^{k_i})(q-\epsilon^{k_i})} 
$$
$$
=\frac{(n-e(1))}{\Phi_h(q)}\frac{\prod_{d|h,d\, prime}p^{\frac{\phi(h)}{p-1}}}{(-1)^{1+\phi(h)/2}h^{\phi(h)}}\det\begin{pmatrix}0 & Q_{\phi(h)}\\M_h^t & C\end{pmatrix},
$$
где у вектора $M_h$ для $j=1,\dots,\phi(h)$ $j$-я координата равна $L_{h,j+1}$.$\blacksquare$ 
\end{prop}

\begin{prop}\label{prop19}
$$
\frac{n}{1-q}=\sum_{d|h} \frac{1/qF_d(q)+q^dF_d(1/q)}{1-q^d}.\blacksquare
$$
\end{prop}

\bigskip

\section{Показатели квазиоднородности простых особенностей}

\bigskip

  В случае $A-D-E$ для показателей квазиоднородности $a, b, c$ связанных с системами корней простых особенностей известна формула \cite{23}
$$
\sum_{i=1}^nq^{e_i}=q^{-h}\frac{(q^h-q^a)(q^h-q^b)(q^h-q^c)}{(q^a-1)(q^b-1)(q^c-1)}.
$$
Формула (1) с учётом соотношений $a+b+c=h+1$ и $c=h/2$ приводит к равенству
$$
B(q)=\frac{1}{q-1}\left(\frac{q(q^{h-a}-1)(q^{h-b}-1)}{(q^a-1)(q^b-1)}-n\right).       \eqno (19)
$$
С помощью ещё одного известного соотношения $2ab=g$ ($g$ -- порядок отвечающей простой особенности бинарной полиэдральной группы) можно заменить $a$ и $b$ в этой формуле:
$$
a,b=\frac{h+2\pm\sqrt{(h+2)^2-8g}}{4}.
$$
Если $a<b$, то $a$ совпадает с размерностью представления бинарной полиэдральной группы, отвечающего согласно соответствию Маккея вершине ветвления диаграммы Дынкина (в случае $A_n$ -- с размерностью любого неприводимого представления). Также можно заменить параметры в формуле (19) с помощью соотношения из \cite{23} для $A_{2n-1}$, $B_n$, $C_n$, $D_n$, $E_6$, $E_7$, $E_8$
$$
v=\frac{h}{abc}=\frac{c}{\alpha\beta\gamma}.
$$
где $v$ -- сферический объём отвечающего системе корней треугольника Шварца, $c$ -- определитель матрицы Картана, $\alpha,\beta,\gamma$ -- длины ветвей диаграммы Дынкина.

\bigskip

\section{Многочлен Дынкина и перечисление антицепей в $\Phi^+$}

\bigskip

  Определим многочлен $D_\Phi(q)$ равенством
$$
D_\Phi(q):=\frac{1-q^h}{1-q}\sum_{i=1}^nq^{e_i-1}.       \eqno (20)
$$
В случае $A-D-E$ $D_\Phi(q)$ появляется в \cite{22} как многочлен Гильберта $W$-действия перестановками на системе корней $\Phi$, в \cite{19} как многочлен Дынкина представления соответствующей группы Ли, в \cite{9} как скалярный квадрат обобщённого характера Холла-Литлвуда и в \cite{9}, \cite{19} как частное функций роста группы Вейля $W$ и её подгруппы. Из формул (1) и (20) следует равенство
$$
B(q)=\frac{qD_\Phi(q)}{q^h-1}-\frac{n}{q-1}. 
$$
В случае $A-D-E$ похожий многочлен появляется в \cite{18}, р. 1206, при перечислении антицепей для частичного порядка на $\Phi^+$
$$
M_\Phi(q):=\frac{1-q^{h-1}}{1-q}\sum_{i=1}^nq^{e_i-1}.
$$
Аналогично имеем равенство
$$
B(q)=\frac{qM_\Phi(q)}{q^{h-1}-1}-\frac{n}{q-1}. 
$$

\end {document}